\documentclass[12pt]{article}

\usepackage[latin1]{inputenc}
\usepackage{amsmath,amssymb}
\usepackage{latexsym}

\newtheorem{theorem}{Theorem}[section]
\newtheorem{corollary}[theorem]{Corollary}
\newtheorem{lemma}[theorem]{Lemma}
\newtheorem{proposition}[theorem]{Proposition}

\newtheorem{remark}[theorem]{Remark}
\newtheorem{example}[theorem]{Example}
\numberwithin{equation}{section}

\def\proof{{\medskip\noindent {\bf Proof. }}}

\def\qed{{\hfill $\square$ \bigskip}}

\def\square{{\vcenter{\vbox{\hrule height.3pt
        \hbox{\vrule width.3pt height5pt \kern5pt
           \vrule width.3pt}
        \hrule height.3pt}}}}

 \def\sE {{\cal E}}

\def\wt{\widetilde}

\def\ol{\overline}
\def\E{{\mathbb E}}
\def\P{{\mathbb P}}
\def\norm#1{{\Vert #1 \Vert}}

\def\lam{{\lambda}}
\def\brack#1{{\langle #1 \rangle}}
\def\bee{\begin{equation}}
\def\eee{\end{equation}}

\def\R{{\mathbb R}}
\def\E{{{\mathbb E}\,}}
\def\P{{\mathbb P}}

\def\Z{{\mathbb Z}}

\def\lam{{\lambda}}

\def\al{{\alpha}}
\def\grad{{\nabla}}
\def\proof{{\medskip\noindent {\bf Proof. }}}

\def\qed{{\hfill $\square$ \bigskip}}

\def\eps{\varepsilon}

\def\angel#1{{\langle#1\rangle}}
\def\norm#1{\Vert #1 \Vert}

 \def\qq {\qquad}

\def\wt{\widetilde}
\def\ol{\overline}

\def\ni{\noindent }
\def\ms{\medskip}

\def\square{{\vcenter{\vbox{\hrule height.3pt
        \hbox{\vrule width.3pt height5pt \kern5pt
           \vrule width.3pt}
        \hrule height.3pt}}}}

\def\tfrac#1#2{{\textstyle {\frac{#1}{#2}}}}

\def\tlint{{- \kern-0.85em \int \kern-0.2em}}  
\def\dlint{{- \kern-1.05em \int \kern-0.4em}}  

 \def\sE {{\cal E}}

\def\nn{{\nonumber}}

\begin{document}

\title{Symmetric jump processes: localization, heat kernels, and convergence}

\author{Richard F. Bass\thanks{Research partially supported by NSF grant
DMS0601783.}, Moritz Kassmann
\thanks{Research partially supported by DFG (German Science
Foundation) through SFB 611.}, and Takashi Kumagai
\thanks{Research partially supported by the Grant-in-Aid for Scientific
Research (B) 18340027 (Japan).} 
}

\maketitle

\begin{abstract}  
\noindent {\it Abstract:} We consider symmetric processes of pure jump type.
We prove local estimates on the probability of exiting balls, the  
H\"older continuity of harmonic functions and of heat kernels, and 
convergence of a sequence of such processes.

\end{abstract}

\section{Introduction}\label{S:I}

Suppose $J:\R^d\times \R^d\to [0,\infty)$ is a symmetric function
satisfying
$$
\frac{c_1}{|y-x|^{\beta_1}}\leq J(x,y)\leq \frac{c_2}{|y-x|^{\beta_2}}
$$
if $|y-x|\leq 1$ and 0 otherwise.
Define the Dirichlet form
\begin{equation}\label{DFdedef}
\sE(f,f)=\int\int (f(y)-f(x))^2 J(x,y)\, dy\, dx,
\end{equation}
and we take as the domain of $\sE$ the closure with respect to the
norm $(\norm{f}_{L^2(\R^d)}$ $+\sE(f,f))^{1/2}$
of the Lipschitz  
functions with compact support. When $\beta_1=\beta_2$, the
Dirichlet form and associated infinitesimal generator are said to be
of fixed order, namely, $\beta_1$, while if $\beta_1<\beta_2$, the
generator is of variable order. The variable order case allows
for considerable variability in the jump intensities and directions.
 
In \cite{BBCK} a number of results were proved for the Hunt process $X$
associated with $\sE$, including exit probabilities, heat kernel estimates,
a parabolic Harnack inequality, and the lack of continuity of harmonic functions. 
The last is perhaps the most interesting: it was shown that there exist
harmonic functions that are not continuous. 

This paper could be considered a sequel to \cite{BBCK}, although the set
of authors for the present paper neither contains nor is contained in the
set of authors of \cite{BBCK}. We prove three main results, which we
discuss in turn.

First we discuss estimates on exit probabilities.  In \cite{BBCK} some
 estimates were obtained on $\P^x(\tau_{B(x,r)}<t)$. These estimates
held for all $x$, but were very crude, and were not sensitive to the behavior 
of $J(x,y)$ when $y$ is close to $x$. We show in the current paper
that to a large extent the behavior of these exit probabilities depend
on the size of $J(x,y)$ for $y$ near $x$. We also allow large jumps,
which translates to allowing $J(x,y)\ne 0$ for $|y-x|>1$. 
In Example \ref{Ex1} we
show how under some smoothness in $J$, we can get fairly precise
estimates.

Our motivation for obtaining better bounds on exit probabilities is
to consider the question of when harmonic functions and the heat kernel
are continuous. The example in \cite{BBCK} shows this continuity need
not always hold. However, when $J$ possesses a minimal amount of smoothness,
we establish that indeed harmonic functions are H\"older continuous,
and the heat kernel is also H\"older continuous. The technique for
showing the H\"older continuity of harmonic functions is based on ideas from
\cite{BKa05}, where the non-symmetric case was considered. More interesting is
the part of the proof where we show that H\"older continuity of harmonic functions
plus global bounds on the heat kernel imply H\"older continuity of the
heat kernel. This argument is of independent interest, and should be applicable
in many other situations.

Finally, we suppose we have a sequence of functions $J_n$ with corresponding
Dirichlet forms and Hunt processes. We show that if the $J_n$ converge
weakly to $J$, and some uniform integrability holds, then the corresponding 
processes converge. Note only weak convergence of the $J_n$ is needed.
This is in contrast to the diffusion case, where it is known that
weak convergence is not sufficient, and a much stronger type of convergence
of the Dirichlet forms is required; see \cite{SZ}. 

Our assumptions and results are stated and proved in the next three
sections, the exit probabilities in Section \ref{S:EP}, the regularity
in Section \ref{S:R}, and the weak convergence in Section \ref{S:WC}.

\section{Exit probabilities}\label{S:EP}

Suppose $J:\R^d\times \R^d\to [0,\infty)$ is jointly measurable. We suppose
throughout this paper that
there exist constants $\kappa_1,\kappa_2, \kappa_3>0$ and $\beta_1, \beta_2\in (0,2)$
such that
\bee\label{E01}
\frac{\kappa_1}{|x-y|^{d+\beta_1}}\leq J(x,y)\leq \frac{\kappa_2}{|x-y|^{d+\beta_2}},
\qq |x-y|\leq 1,
\eee
and 
\bee\label{E02}
\int_{|x-y|>1} J(x,y)\, dy\leq \kappa_3, \qq x\in \R^d.
\eee
The constants $\beta_1, \beta_2, c_1, c_2, c_3$ play only a limited
role in what follows and \eqref{E01} and \eqref{E02} are used
to guarantee a certain amount of regularity. Much more important is the
$\al$ that is introduced in \eqref{E1}.
Define a Dirichlet form $\sE=\sE_J$ by
\bee\label{E03}
\sE(f,f)=\int\int (f(y)-f(x))^2 J(x,y)\, dy\, dx,
\eee
where we take the domain to be the closure of the 
Lipschitz functions with 
compact support with respect to the norm $(\norm{f}_2+(\sE(f,f))^{1/2}$.
Let $X$ be the Hunt process associated with the Dirichlet form $\sE$.
Let $B(x,r)$ denote the open ball of radius $r$ centered at $x$. 
The letter $c$ with or without subscripts will denotes constants whose
exact values are unimportant and which may change 
from line to line.

We remark that if we define $J_1(x,y)=J(x,y)1_{(|x-y|\leq 1)}$ and
define the corresponding Dirichlet form in terms of $J_1$, then
the Hunt process $X^{(1)}$ corresponding to this Dirichlet form is conservative
by \cite[Theorem 1.1]{BBCK}. 
Using a construction due to Meyer
(see \cite[Remark 3.4]{BBCK} and \cite[Section 3.1]{BGK}) 
we can use $X^{(1)}$
 to obtain $X$. This is a probabilistic procedure that involves
adding jumps. Only finitely many jumps are added in any finite time
interval, and we deduce from this 
construction that $X$ is also conservative.

We now fix $z_0\in \R^d$  and
assume that there exist  constants  $\kappa_4$ and 
$\al\in (0,2)$  such that
\bee\label{E1}
J(x,y) \geq \kappa_4 |x-y|^{-d-\al}, \qq x,y\in B(z_0,3r).
\eee
Here $\al$ may depend on $z_0$.

Define
\bee\label{E2}
L_1(x,s)=\int_{|x-w|\geq s} J(x,w)\, dw,
\eee
\bee\label{E3}
L_2(x,s) =\int_{|x-w|\leq s} |x-w|^2  J(x,w)\, dw,
\eee
and let 
\bee\label{E4}
L(z_0,r)
=\sup_{x\in B(z_0,3r)} L_1(x,r)+\sup_{x\in B(z_0,3r)} \sup_{s\leq r}s^{d} [s^{-2}L_2(x,s)]^{\frac{d+\al}{\al}}.
\eee
From \eqref{E1} we see that
\bee\label{E4.5}
L(z_0,r)\geq cr^{-\al}.
\eee

\begin{theorem}\label{T1}
Suppose \eqref{E01}, \eqref{E02}, and \eqref{E1} hold.
There exists $c_1$ (depending only on $d$, $\kappa_4$, and $\al$)  
such that if $r\in (0,1)$, then 
for $x\in B(z_0,r)$,
$$\P^x(\tau_{B(x,r)}<t)\leq c_1tL(z_0,r).$$
\end{theorem}

\proof 
Let $x_0,y_0$ be fixed, let $R=|y_0-x_0|$, and
suppose $R\geq 18(d+\al)r/\al$. 
By \eqref{E4.5}, the result is immediate if $t>r^\al$, so let us
suppose $t\leq r^\al$. 
Define
\bee\label{E5}
\wt J(x,y)=\begin{cases} J(x,y) 
& \mbox{ if }x,y \mbox{ are both in }B(z_0,3r) \mbox{ and } |x-y|< R,\\ 
\kappa_4|x-y|^{-d-\al}& \mbox{ if at least one of  $x$ and $y$ is not in }\\
&\qq B(z_0,3r) \mbox{ and } |x-y|< R,\\
0 &\mbox{ otherwise.}
\end{cases}
\eee
Define
$$\wt J_\delta(x,y)=\wt J(x,y) 1_{(|x-y|\leq \delta)},$$
where we will choose 
$\delta \in [6r, R)$ in a moment.
Let $\wt X$ be the Hunt process corresponding to $\wt J$ and $\wt X^{(\delta)}$
the Hunt process associated with $\wt J_\delta$.

We have the Nash inequality
(see, e.g., (3.9) of \cite{BBCK}):
\bee\label{E6}
\norm{u}_2^{2+\frac{2\al}{d}} \leq 
c\Big(\int\int_{|x-y|<\delta} \frac{(u(x)-u(y))^2}{|x-y|^{d+\al}}\, dy\, dx
+\delta^{-\al}\|u\|_2^2\Big) 
\, \norm{u}_1^{2\al/d}.
\eee
Using \eqref{E1}
we obtain from this that 
\bee\label{E7}
\norm{u}_2^{2+\frac{2\al}{d}} \leq c\Big(\int\int (u(x)-u(y))^2 \wt J_\delta (x,y)\, dy\, dx
+\delta^{-\al} \norm{u}_2^2\Big) \, \norm{u}_1^{2\al/d}.
\eee

Let
\begin{align}  
\delta&=\frac{R\al}{3(d+\al)},\label{E71}\\ 
N(\delta)&=\delta^{-\al}+\sup_{x\in B(z_0,3r)}\delta^{-2}L_2(x,\delta), \label{E72}\\
\lam&=\frac{1}{3\delta} \log(1/(N(\delta)t)). \label{E73}
\end{align}
Let   
$\psi(x)=\lam(R-|x-x_0|)^+$.
Set
$$\Gamma(f,f)(x)=\int (f(y)-f(x))^2 \wt J_\delta(x,y)\, dy.$$
Since $|e^t-1|^2\leq t^2 e^{2t}$, $|\psi(x)-\psi(y)|\leq \lam |x-y|$, and
$\wt J_\delta(x,y)=0$ unless $|x-y|<\delta$, then
\begin{align}
e^{-2\psi(x)}\Gamma(e^\psi, e^\psi)(x)&=\int_{|x-y|\leq \delta} 
\Big( e^{\psi(x)-\psi(y)}-1\Big)^2 \wt J_\delta(x,y)\, dy\nn\\
&\leq e^{2\lam \delta} \lam^2 \int_{|x-y|\leq \delta}
|x-y|^2 \wt J_\delta(x,y)\, dy.\label{E74}
\end{align}

Since  $\delta \geq 6r$, then by our definition of $\wt J$ we have that
the integral on the last line of \eqref{E74}  is bounded by
$\sup_{x\in B(z_0,3r)}L_2(x,\delta)+\delta^{2-\al}.$
We therefore have
\begin{align*}
e^{-2\psi(x)}\Gamma(e^\psi, e^\psi)(x)
&\leq e^{2\lam \delta} \lam^2 \delta^2 N(\delta)\\
&\leq e^{3\lam \delta} N(\delta).
\end{align*}
We obtain in the same way the same upper bound for  $e^{2\psi(x)} \Gamma (e^{-\psi}, e^{-\psi})(x)$. So by \cite[Theorem 3.25]{CKS} we have
\bee\label{E9}
p_\delta(t,x_0,y_0)\leq ct^{-d/\al} e^{ct\delta^{-\al}} e^{-\lam R+c e^{3\lam \delta}N(\delta)t},
\eee
where $p_\delta$ is the transition density for $\wt X^{(\delta)}$.
(Note that by \cite[Theorem 3.1]{BBCK}, the transition density $p_\delta(t,x,y)$ 
exists for $x,y\in \R^d\setminus {\cal N}$, where $\cal N$ is a set of 
capacity zero, called a properly exceptional set. We will take $x_0,y_0\in 
\R^d\setminus {\cal N}$.)    
Since   $t\leq  r^\al\leq c\delta^{\al}$, 
we then get
$$p_\delta(t,x_0,y_0)\leq ct^{-d/\al} e^{-\lam R} =ct^{-d/\al} (N(\delta)t)^{R/3\delta}.$$
Our bound now becomes
\begin{align*}
p_\delta(t,x_0,y_0)&\leq ct^{-d/\al} t^{1+\frac{d}{\al}} 
N(\delta)^{(d+\al)/\al}\\ &= ct N(\delta)^{(d+\al)/\al}.
\end{align*}

Since $\delta\geq 6r$,           
then
$$\norm{\wt J-\wt J_\delta}_\infty 
\leq c\delta^{-(d+\al)},$$
so by \cite[Lemma 3.1]{BGK} and by (\ref{E71}) 
\begin{align*}
p(t,x_0,y_0)&\leq p_\delta(t,x_0,y_0)+ct\delta^{-(d+\al)}\\
&\leq ct[\sup_{x\in B(z_0,3r)}\delta^{-2}
L_2(x,\delta)+\delta^{-\al}]^{\frac{d+\al}{\al}}+ct R^{-(d+\al)}.
\end{align*}
Since
\begin{align*}
\sup_{x\in B(z_0,3r)}\delta^{-2}L_2(x,\delta)+\delta^{-\al}
&\leq c\delta^{-2}[ \sup_{x\in B(z_0,3r)}L_2(x,r)+\delta^{2-\al}]\\
&\leq cR^{-2}\sup_{x\in B(z_0,3r)}L_2(x,r)+cR^{-\al},
\end{align*}
then,
because $\wt X$ is conservative, 
integrating over 
$R\geq r/2$ gives us  
$$\P^{x_0}(|\wt X_t-x_0|\geq r/2)\leq ctr^{d} \Big[r^{-2}
\sup_{x\in B(z_0,3r)}L_2(x,r)\Big]^{\frac{d+\al}{\al}}+ctr^{-\al}\le ctL(z_0,r).$$
By \cite[Lemma 3.8]{BBCK} we then have
$$\P^{x_0}(\sup_{s\leq t} |\wt X_s-x_0|>r)\leq ctL(z_0,r).$$

We now use Meyer's construction to compare  
$\wt X$ to $X$. Using this construction we obtain, 
for $x\in B(z_0,r)$, 
\begin{align*}
\P^x( X_s\ne \wt X_s\mbox{ for some }s\leq t)&\leq 
t\sup_{x'\in B(z_0,2r)}\int_{B(z_0,3r)^c}|J(x',y)-\wt J(x',y)|dy\\
&\le ctL(z_0,r).
\end{align*}
(The first inequality can be obtained by observing the processes 
$X$ and $\wt X$ killed on exiting $B(z_0, 2r)$.)  Therefore, for $x\in B(z_0,r)$,
\begin{align*}
\P^x(\sup_{s\leq t}|X_s-x|>r)
&\leq \P^x(\sup_{s\leq t} |\wt X_s-x|>r)+\P^x(X_s\ne \wt X_s
\mbox{ for some }s\leq t)\\
&\leq ctL(z_0,r).
\end{align*}
\qed

\begin{corollary}\label{C1}
Suppose \eqref{E01} and \eqref{E02} hold.
Suppose instead of \eqref{E1} we have 
\bee\label{E101}
J(x,y)\geq \kappa_4 |x-y|^{-d-\al}-K(x,y), \qq x,y\in B(z_0,3r),
\eee
where 
\bee\label{E102}
\int_{|x-y|\leq \delta} K(x,y)\, dy\leq \kappa_5 \delta^{-\al}
\eee 
for all $x\in B(z_0,3r)$ and all $\delta\leq r$.
Then the conclusion of Theorem \ref{T1} still holds.
\end{corollary}

\proof 
The only place the lower bound on $J(x,y)$ plays a role
is in deriving \eqref{E7} from \eqref{E6}. If we have \eqref{E101} instead of
\eqref{E1}, then in place of \eqref{E7}
we now have
\begin{align}
\norm{u}_2^{2+\frac{2\al}{d}} &\leq c\Big(\int\int (u(x)-u(y))^2 \wt J_\delta (x,y)\, dy\, dx\label{E103}\\
&\qq +\int\int_{|x-y|\leq \delta} (u(x)-u(y))^2 K(x,y)\, dy\, dx
+\delta^{-\al} \norm{u}_2^2\Big) \, \norm{u}_1^{2\al/d}.\nn
\end{align}
But by our assumption on $K(x,y)$, the double integral with $K$ in the integrand is
bounded by
$$ \int \Big(\int_{|x-y|\leq \delta} K(x,y)\, dy\Big)\ u(x)^2\, dx
\leq c\delta^{-\al}\norm{u}_2^2.$$
\qed

\begin{example}\label{Ex1}{\rm 
Suppose $\eps>0$ and there exists a function $s:\R^d\to (\eps, 2-\eps)$
such that 
\bee\label{E1031}
|s(x)-s(y)|\leq c \log(2/|x-y|), \qq |x-y|<1.
\eee
Suppose there exist constants $c_1, c_2$ such that
\bee\label{E1032}\frac{c_1}{|x-y|^{d+s(x)\land s(y)}}\leq J(x,y)
\leq \frac{c_2}{|x-y|^{d+s(x)\lor s(y)}}.
\eee
Suppose further that (\ref{E02}) holds. 
We show that $L(z_0,r)$ is comparable to $r^{-s(z_0)}$
if $r<1$. 

To see this, note that
\bee\label{E104}
|x-y|^{s(x)-s(y)}\leq 
|x-y|^{-c/\log(2/|x-y|)}\leq e^c,
\eee
if $|x-y|\leq 1$
and similarly we have 
\bee\label{E105}|x-y|^{s(x)-s(y)}\geq 
|x-y|^{c/\log(2/|x-y|)}\geq e^{-c}.
\eee
If we fix $x$ and let
$$M(v)=\sup_{|x-w|=v} J(x,w),$$
then for $v\leq 1$ 
\begin{align*}
M(v)&\leq \sup_{|x-w|=v} \frac{c}{v^{d+s(x)}} v^{-|s(x)-s(w)|}\\
&\leq \frac{c}{v^{d+s(x)}}.
\end{align*}
We then estimate for $r\leq 1$
\begin{align*}
L_2(x,r)&\leq c\int_0^r v^2 M(v) v^{d-1}\, dv\\
&\leq c\int_0^r v^{1-s(x)}\, dv=cr^{2-s(x)},
\end{align*}
We can similarly obtain a bound for $L_1(x,r)$:
\begin{align*}
L_1(x,r)&\leq c\int_r^1 M(v)v^{d-1}\, dv+\int_{|x-w|>1}J(x,w)\, dw\\
&\leq c\int_r^1 v^{-1-s(x)} \, dv+c\\
&\leq cr^{-s(x)}+c\leq cr^{-s(x)}
\end{align*}
if $r\leq 1$.

Next, for $x\in B(z_0,3r)$, we have $r^{-s(x)}$ is
comparable to $r^{-s(z_0)}$ for $r\leq 1$. 
To see this,
$$c\leq r^{s(x)-s(z_0)}\leq 
r^{-|s(x)-s(z_0)|}\leq c'$$
as in \eqref{E104} and \eqref{E105}.

If we take $\al$ in \eqref{E1} to be $\inf_{x\in B(z_0,3r)} s(x)$, then we
conclude
$$L(x,r)\leq cr^{-s(z_0)}+cr^{-\al}\leq cr^{-s(z_0)},$$
so $$\P^x(\tau_r\leq t)\leq ct r^{-s(z_0)}, \qq x\in B(z_0,r).$$
}
\end{example}

\section{Regularity}\label{S:R}

We suppose throughout this section that \eqref{E01} and \eqref{E02} hold.
We suppose in addition first  that there exists $c$ such that
\bee\label{C01} \int_A J(z,y)\, dy\geq c L(x,r)
\eee
whenever $r\in (0,1)$, 
$A\subset B(x,3r)$, $|A|\geq \frac13 |B(x,r)|$, 
$x\in \R^d$, and $z\in B(x,r/2)$
and second there exist $\sigma$ and $c$ such that
\bee\label{C02}
\frac{L_1(x,\lam r)}{L_1(x,r)}\leq c\lam^{-\sigma}, \qq x\in \R^d, r\in (0,1),
\lam\in (1,1/r).
\eee
It is easy to check that \eqref{C01} and \eqref{C02} hold for Example \ref{Ex1}.

We say a function $h$ is harmonic in a ball $B(x_0,r)$ if 
$h(X_{t\land \tau_{B(x_0,r(1-\eps))}})$ is a $\P^x$ martingale for q.e. 
$x$ and every $\eps\in (0,1)$.

\begin{theorem}\label{T.R1}
Suppose \eqref{E01}, \eqref{E02}, \eqref{C01}, and \eqref{C02} hold.
There exist $c_1$ and $\gamma$ such that if 
$h$ is bounded in $\R^d$ and harmonic in
a ball $B(x_0,r)$, then
\bee\label{C03}
|h(x)-h(y)|\leq c_1\Big(\frac{|x-y|}{r}\Big)^\gamma\norm{h}_\infty, \qq x,y\in B(x_0,r/2).
\eee
\end{theorem}

\proof As in \cite{CK03, CK07} we have the 
L\'evy system formula:
\bee\label{C04}
\E^x\Big[\sum_{s\leq T} f(X_{s-},X_s)\Big]=
\E^x\Big[\int_0^T \Big(\int f(X_s,y) J(X_s,y)\, dy\Big)\, ds\Big]
\eee
for any nonnegative $f$ that is 0 on the diagonal, for every bounded stopping
time $T$, and q.e. starting point $x$.
Given this, the proof is nearly identical to that in 
\cite[Theorem 2.2]{BKa05}.
\qed

We obtain a crude estimate on the expectation of the exit times.

\begin{lemma}\label{L.R2}
Assume the lower bound of  \eqref{E01}. Then 
there exists $c_1$ such that
$$\E^x \tau_r\leq c_1 r^{\beta_1}, \qq x\in \R^d, r\in (0,1/2).$$
\end{lemma}

\proof
The expression $\sum_{s\leq t\land \tau_r} 1_{(|X_s-X_{s-}|>2r)}$ is
1 if there is a jump of size at least $2r$ before time $t\land \tau_r$,
in which case the process exits $B(x,r)$ before or at time $t$, or 0 if there is no such jump. 
So
\begin{align*}
\P^x(\tau_r\leq t)&\geq \E^x \sum_{s\leq t\land \tau_r} 1_{(|X_s-X_{s-}|>2r)}\\
&=\E^x \int_0^{t\land \tau_r} \int_{B(x,2r)^c} J(X_s,y)\, dy\, ds\\
&\geq cr^{-\beta_1} \E^x[t\land \tau_r]\\
&\geq cr^{-\beta_1}t\P^x(\tau_r>t),
\end{align*}
using the lower bound of \eqref{E01}. 
Thus
$$\P^x(\tau_r>t)\leq 1-cr^{-\beta_1}t \P^x(\tau_r>t),$$
or
$\P^x(\tau_r>t)\leq 1/2$ if we take $t=c^{-1}r^{\beta_1}$.
This holds for every $x\in \R^d$.
Using the Markov property at time $mt$,
$$\P^x(\tau_r>(m+1)t)\leq \E^x[\P^{X_{mt}}(\tau_r>t); \tau_r>mt]\leq \tfrac12 \P^x(\tau_r>mt).$$
By induction  $\P^x(\tau_r>mt)\leq 2^{-m}.$ With this choice of $t$, our lemma follows.
\qed

We next show $\lam$-potentials are H\"older continuous.
Let
$$U^\lam f(x)=\E^x\int_0^\infty e^{-\lam t} f(X_t)\, dt.$$

\begin{proposition}\label{PR3}
Under the same assumption as in Theorem \ref{T.R1}, 
there exist $c_1=c_1(\lambda)$  
and $\gamma'$ such that if $f$ is bounded, then 
$$|U^\lam f(x)-U^\lam f(y)|\leq c_1 |x-y|^{\gamma'} \norm{f}_\infty.$$
\end{proposition}

\proof 
Fix $x_0$, let $r\in (0,1/2)$, and suppose $x,y\in B(x_0,r/2)$.
By the strong Markov property,
\begin{align*}
U^\lam f(x)&=\E^x\int_0^{\tau_r} e^{-\lam t} f(X_t) \, dt
+\E^x (e^{-\lam \tau_r}-1)U^\lam f(X_{\tau_r})\\
&\qq + \E^x U^\lam f(X_{\tau_r})\\
&=I_1+I_2+I_3,
\end{align*}
and similarly when $x$ is replaced by $y$. We have by Lemma \ref{L.R2}
$$|I_1|\leq \norm{f}_\infty \E^x \tau_r\leq cr^{\beta_1}\norm{f}_\infty$$
and by the mean value theorem and Lemma \ref{L.R2}
$$|I_2|\leq \lam \E^x \tau_r \norm{U^\lam f}_\infty\leq cr^{\beta_1}\norm{f}_\infty,$$ and similarly when $x$ is replaced by $y$.
So
\bee\label{C05}
|U^\lam f(x)-U^\lam f(y)|\leq cr^{\beta_1}\norm{f}_\infty+|\E^x U^\lam f(X_{\tau_r})
-\E^y U^\lam f(X_{\tau_r})|.
\eee
But $z\to \E^z U^\lam f(X_{\tau_r})$ is 
bounded in $\R^d$ and harmonic in $B(x_0,r)$, so 
by Theorem \ref{T.R1}
the second 
term in \eqref{C05} is bounded by
$$c\Big(\frac{|x-y|}{r}\Big)^\gamma \norm{U^\lam f}_\infty.$$
If we use $\norm{U^\lam f}_\infty\leq \frac1{\lam}\norm{f}_\infty$ and
set $r=|x-y|^{1/2}$, then 
\bee\label{C06}|U^\lam f(x)-U^\lam f(y)|\leq (c|x-y|^{\beta_1/2}+c|x-y|^{\gamma/2})\norm{f}_\infty,
\eee
and our result follows.
\qed

Using the spectral theorem, there exists projection operators $E_\mu$ on  the space
$L^2(\R^d, dx)$
such that 
\begin{align} 
f&=\int_0^\infty \, dE_\mu(f),\nn\\
 P_tf&=\int_0^\infty e^{-\mu t} \, dE_\mu(f),\nn\\
 U^\lam f&=\int_0^\infty \frac{1}{\lam+\mu} \, dE_\mu(f).\label{C07}
\end{align}

\begin{proposition}\label{PR4}
Under the same assumptions as in Theorem \ref{T.R1}, 
if $f$ is in $L^2$, then $P_tf$ is equal a.e. to a function that is H\"older
continuous.
\end{proposition}

\proof Write $\brack{f,g}$ for the inner product in $L^2$. Note that in 
what follows $t$ is fixed. Each of our constants may depend on $t$. If $X^{(1)}$ is the Hunt process
associated with the Dirichlet form defined in terms of the kernel
$J_1(x,y)=J(x,y)1_{(|x-y|<1)}$, we  know from 
\cite[Theorem 2.1]{BBCK} 
that $X^{(1)}$ has a transition density $p(t,x,y)$ bounded by $c$. 
Using \cite[Lemma 3.1]{BGK} and Meyer's construction (cf.\ \cite[Section 3]{BBCK}), 
we then can conclude that $X$ also has a transition density bounded by
$c$.
Define
$$h=\int_0^\infty (\lam+\mu)e^{-\mu t}\, dE_\mu(f).$$
Since
$\sup_\mu (\lam+\mu)^2 e^{-2\mu t}\leq c$, then 
$$\int_0^\infty (\lam+\mu)^2 e^{-2\mu t} \, d\brack{E_\mu(f), E_\mu(f)}
\leq c\int_0^\infty \, d\brack{E_\mu(f), E_\mu(f)}=c\norm{f}_2^2,$$
we see that $h$ is a well defined function in $L^2$.

Suppose $g\in L^1$. Then $\norm{P_t g}_1\leq \norm{g}_1$ by Jensen's inequality,
and
$$|P_tg(x)|=\Big|\int p(t,x,y)g(y)\, dy\Big|\leq c\norm{g}_1$$
by the fact that $p(t,x,y)$ is bounded.
So $\norm{P_tg}_\infty\leq c\norm{g}_1$, and it follows that 
$\norm{P_tg}_2\leq c\norm{g}_1$.
Using Cauchy-Schwarz and the fact that
$$\sup_\mu (\lam+\mu) e^{-\mu t/2}\leq c<\infty,$$
we have
\begin{align*}
\brack{h,g}&=\int_0^\infty (\lam+\mu)e^{-\mu t}\, d\brack{E_\mu(f), E_\mu(g)}\\
&\leq \Big(\int_0^\infty (\lam+\mu)e^{-\mu t} \, d\brack{E_\mu(f), E_\mu(f)}\Big)^{1/2}\\
&\qq\qq\times
\Big(\int_0^\infty (\lam+\mu)e^{-\mu t} \, d\brack{E_\mu(g), E_\mu(g)}\Big
)^{1/2}\\
&\leq c\Big(\int_0^\infty  \, d\brack{E_\mu(f), E_\mu(f)}\Big
)^{1/2}
\Big(\int_0^\infty e^{-\mu t/2} \, d\brack{E_\mu(g), E_\mu(g)}\Big
)^{1/2}\\
&=c \norm{f}_2\norm{P_{t/2}g}_2\\
&\leq c \norm{f}_2 \norm{g}_1.
\end{align*}
Taking the supremum over $g\in L^1$ with $L^1$ norm less than 1, 
$\norm{h}_\infty\leq c\norm{f}_2$. But by 
\eqref{C07}
$$U^\lam h=\int_0^\infty e^{-\mu t} \, dE_\mu(f)=P_tf,\qq a.e.,$$
and the H\"older continuity of $P_tf$ follows by Proposition \ref{PR3}.
\qed

Finally we have

\begin{theorem}\label{PT6}
Under the same assumption as in Theorem \ref{T.R1}, 
we can choose $p(t,x,y)$ to be jointly continuous.
\end{theorem}

\proof Fix $y$ and let $f(z)=p(t/2,z,y)$. $f$ is bounded by $c$ (depending on $t$)
and has $L^1$ norm equal to 1, hence $f\in L^2$ with norm bounded by $c$.
Note 
$$P_{t/2}f(x)=\int p(t/2,x,z) f(z)\, dz=\int p(t/2,x,z)p(t/2,z,y)
=p(t,x,y).$$
Using Proposition \ref{PR4} shows that $p(t,x,y)$ is H\"older
continuous with constants independent of $x$ and $y$. This and symmetry gives the result.
\qed

\begin{remark}\label{R3.6}{\rm
The argument we gave deriving the H\"older continuity of the transition
densities from the boundedness of the transition densities plus the
H\"older continuity of harmonic functions holds in much more general
contexts than just jump processes in $\R^d$.}
\end{remark}

\section{Convergence}\label{S:WC}

Suppose now that we have a sequence of jump kernels $J^n(x,y)$
satisfying \eqref{E01}, \eqref{E02}, \eqref{C01}, and \eqref{C02} with
constants independent of $n$. Suppose in addition that
\bee\label{WCA1}
\limsup_{\eta\to 0} \sup_{n,x} \int_{|y-x|\geq \eta^{-1}} J_n(x,y)\, dy\, dx=0,
\eee
\bee\label{WCA2}
\limsup_{\eta\to 0} \sup_{n,x} \int_{|y-x|\leq \eta} |y-x|^2 J_n(x,y)\, dy\, dx=0,
\eee
and
for almost every $\eta$
\bee\label{WCA3}
J_n(x,y)1_{(\eta,\eta^{-1})}(|y-x|)\, dx\, dy\to
J(x,y)1_{(\eta,\eta^{-1})}(|y-x|)\, dx\, dy
\eee
weakly as $n\to \infty$.

Let $\sE^n$ be the Dirichlet forms defined in terms of the $J^n$ with
$P_t^n$,  $U^\lam_n$, and $\P^x_n$  the associated semigroup, resolvent, and
probabilities.
Let $P_t$, $U^\lam$, and $\P^x$ be the semigroup, resolvent, and probabilities
corresponding to the Dirichlet form $\sE_J$ defined in terms of the kernel
$J$.

Under the above set-up we have

\begin{theorem}\label{WT1}
If $f$ is bounded and continuous, 
then $P_t^n f$ converges uniformly on compacts to
$P_tf$. For each $t$,  for  q.e.\ $x$, $\P^x_n$ converges weakly to $\P^x$
with respect to the space $D([0,t])$.
\end{theorem}

\proof The first step is to show that any subsequence $\{n_j\}$ has a further 
subsequence $\{n_{j_k}\}$  such that $U^\lam_{n_{j_k}}f$ converges uniformly
on compacts whenever $f$ is bounded and continuous. The proof
of this is very similar to that of \cite[Proposition 6.2]{BKu08}, and
we refer the reader to that paper.

Now suppose we have a subsequence $\{n'\}$ such that the $U^\lam_{n'} f$ are 
equicontinuous and converge uniformly on compacts whenever $f$ is bounded and
continuous with compact support. Fix such an $f$ and let $H$ be the limit
of $U^\lam_{n'}f$. We will show
\bee\label{WCE1}
\sE_J(H,g)=\angel{f,g}-\lam\angel{H,g}
\eee
whenever $g$ is a Lipschitz function  
with compact support, where $\sE_J$ is the Dirichlet
form corresponding to the kernel $J$. This will prove that $H$ is the 
$\lam$-resolvent of $f$ with respect to $\sE_J$, that is, $H=U^\lam f$. 
We can then conclude that the full sequence $U^\lam_n f$ converges to
$U^\lam f$ whenever $f$ is bounded and continuous with compact support. 
The assertions about the convergence of $P_t^n$ and $\P^x_n$ then follow
as in \cite[Proposition 6.2]{BKu08}.

So we need to prove $H$ satisfies \eqref{WCE1}. We drop the primes for
legibility.

We know
\bee\label{WCE3}
\sE^n(U^\lam_n f, U^\lam_n f)=\angel{f,U^\lam_n f}-\lam \angel{U^\lam_n f,
U^\lam_n f}.
\eee
Since $\norm{U^\lam_n f}_2\leq (1/\lam) \norm{f}_2$ (by Jensen's inequality),
we have by Cauchy-Schwarz that
$$\sup_n \sE^n (U^\lam_n f, U^\lam_n f)\leq c<\infty.$$
Since the $U^\lam_n f$ are equicontinuous and converge uniformly to
$H$ on $B(0, \eta^{-1})-\ol{B(0,\eta)}$ for almost every $\eta$, then
\begin{align*}
\int\int_{\eta<|y-x|<\eta^{-1}} &(H(y)-H(x))^2 J(x,y)\, dy\, dx\\
&\leq \limsup_{n\to \infty} \int\int_{\eta<|y-x|<\eta^{-1}} 
(U^\lam_n f(y)-U^\lam_n f(x))^2 J_n(x,y)\, dy\, dx\\
&\leq \limsup_n \sE^n(U^\lam_n f, U^\lam_n f)\leq c<\infty.
\end{align*}
Letting $\eta\to 0$ (while avoiding the null set), we have
\bee\label{WCE2}
\sE_J(H,H)<\infty.
\eee

Fix a Lipschitz function $g$ 
with compact support and choose $M$ large enough so that
the support of $g$ is contained in $B(0,M)$. Then
\begin{align*}
\Big| \int\int_{|y-x|\geq \eta^{-1}} &
(U^\lam_n f(y)-U^\lam_n f(x))(g(y)-g(x)) J_n(x,y) \, dy\, dx\Big|\\
&\leq \Big( \int\int (U^\lam_n f(y)-U^\lam_n f(x))^2 J_n(x,y)\, dy\, dx\Big)^{1/2}\\
&\qq\times \Big(\int\int_{|y-x|\geq \eta^{-1}} (g(y)-g(x))^2 J_n(x,y)\, dy\, dx\Big)^{1/2}.
\end{align*}
The first factor is $(\sE^n(U^\lam_n f, U^\lam_nf))^{1/2}$, while the
second factor is bounded by
$$\norm{g}_\infty \Big(\int_{B(0,M)}\int_{|y-x|\geq \eta^{-1}} J_n(x,y)\, dx\, dy\Big)^{1/2},$$
which, in view of \eqref{WCA1}, will be small if $\eta$ is small. 
Similarly,
\begin{align*}
\Big| \int\int_{|y-x|\leq \eta} &
(U^\lam_n f(y)-U^\lam_n f(x))(g(y)-g(x)) J_n(x,y) \, dy\, dx\Big|\\
&\leq \Big( \int\int (U^\lam_n f(y)-U^\lam_n f(x))^2 J_n(x,y)\, dy\, dx\Big)^{1/2}\\
&\qq\times \Big(\int\int_{|y-x|\leq \eta} (g(y)-g(x))^2 J_n(x,y)\, dy\, dx\Big)^{1/2}.
\end{align*}
The first factor is as before, while the second is bounded by
$$\norm{\grad g}_\infty \Big(\int_{B(0,M)} \int_{|y-x|\leq \eta}|y-x|^2
J_n(x,y)\, dx\, dy\Big)^{1/2}.$$
In view of \eqref{WCA2}, the second factor will be small if $\eta$ is small.
Similarly, using \eqref{WCE2}, we have
$$\Big|\int\int_{|y-x|\notin (\eta,\eta^{-1})}
(H(y)-H(x))(g(y)-g(x)) J(x,y)\, dy\, dx\Big|$$
will be small if $\eta$ is taken small enough. 

By \eqref{WCA3}, \eqref{E01}, \eqref{E02}, and the fact that the $U^\lam_nf$ are
equicontinuous and converge to $H$ uniformly on compacts,
for almost every $\eta$
\begin{align*}
\int\int_{|y-x|\in (\eta, \eta^{-1})}& (U^\lam_nf(y)-U^\lam_nf(x))(g(y)-g(x))J_n(x,y)\, dy\, dx\\
&\to \int\int_{|y-x|\in (\eta, \eta^{-1})} (H(y)-H(x))(g(y)-g(x))J(x,y)\, dy\, dx.
\end{align*}
It follows that
\bee\label{WCE4}
\sE^n(U^\lam_nf,g)\to \sE_J(H,g).
\eee
But
$$\sE^n(U^\lam_nf,g)=\angel{f,g}-\lam\angel{U^\lam_nf,g}\to 
\angel{f,g}-\lam\angel{H,g}.$$
Combining with \eqref{WCE4} proves \eqref{WCE1}.
\qed

\begin{remark}\label{WR1} {\rm 
One can modify the above proof to obtain  
a central limit theorem for symmetric Markov chains.
Suppose for each $n$ we have a symmetric Markov chain on $n^{-1}\Z^d$ with unbounded
range with
conductances $C^n_{xy}$. If $\nu_n$ is the measure that gives
mass $n^{-d}$ to each point in $n^{-1}\Z^d$, then we can define
the Dirichlet form
$$\sE_n(f,f)=\sum_{x,y\in n^{-1}\Z^d} (f(x)-f(y))^2 C^n_{xy}$$
with respect to the measure $\nu_n$. 
Under appropriate assumptions
analogous to those in 
Sections \ref{S:EP} and \ref{S:R}, 
one can show that the semigroups corresponding to
$\sE_n$ converge to those of  $\sE$
and in addition there is weak convergence of the probability laws.
Since the details are rather lengthy, we leave this 
to the interested reader. 
}
\end{remark} 
\begin{remark}\label{WR121} {\rm 
We can also prove the following approximation of a jump process 
by Markov chains, which is a generalization of \cite[Theorem 2.3]{HK07}. 
 Suppose $J:\R^d\times \R^d\to [0,\infty)$ is a symmetric measurable 
function satisfying \eqref{E01}, \eqref{E02}, \eqref{C01}, and \eqref{C02}.  
Define the conductivity functions $C^n\colon n^{-1}\Z^d \times n^{-1}\Z^d \to[0,\infty)$ by
\[
C^n(x,y) = 
n^{2d}\int_{|x-\xi|_\infty<  \frac 1{2n}}\int_{|y-\zeta|_\infty< \frac 1{2n}} 
J(\xi,\zeta) d\xi\,d\zeta ~\text{ for }~ x\ne y\in n^{-1}\Z^d,
\]
and $C^n(x,x) = 0$, where $|x-y|_\infty=\max_{1\le i\le d}|x_i-y_i|$. 
Let $X$ be the Hunt process corresponding to the Dirichlet form 
given by (\ref{DFdedef}).
Then the sequence of processes corresponding to $C^n$ converges weakly to $X$. 
Given Remark \ref{WR1}, the proof is standard. 
}
\end{remark} 

\begin{remark}\label{WR221} {\rm 
As we mentioned at the beginning of Section 2, the assumptions 
\eqref{E01} and \eqref{E02} are used 
to guarantee a certain amount of regularity, namely, conservativeness and 
the existence of the heat kernel. However, one can relax these assumptions. 
All of the  results in this paper hold if instead of \eqref{E01} 
and \eqref{E02} we assume \eqref{E02}, 
\eqref{E1} for all $z_0\in \R^d$ and the following:
\[
\int_{|x-y|\le 1} |x-y|^2J(x,y)\, dy\leq \kappa_5, \qq x\in \R^d.
\]
}
\end{remark}


\medskip

\ni {\bf Richard F. Bass}\\
Department of Mathematics\\
University of Connecticut \\
Storrs, CT 06269-3009, USA\\
{\it bass@math.uconn.edu}
\ms

\ni {\bf Moritz Kassmann}\\
Institut f\"{u}r Angewandte Mathematik\\
Universit\"{a}t Bonn \\
Beringstrasse 6\\
D-53115 Bonn, Germany\\
{\it kassmann@iam.uni-bonn.de}
\ms

\ni {\bf Takashi Kumagai}\\
Department of Mathematics\\
Kyoto University\\
Kyoto 606-8502, Japan\\
{\it  kumagai@math.kyoto-u.ac.jp}

\end{document}